\documentclass[11pt,a4paper,twoside]{article}
\usepackage[T1]{fontenc}
\usepackage{amssymb}
\usepackage{amsfonts}
\usepackage{amsmath,amsthm}
\usepackage{fancyhdr}
\usepackage{url}
\usepackage{booktabs}
\usepackage{float}
\theoremstyle{definition}

\begin{document}
\begin{center}
\vspace*{2pt}
{\Large \textbf{A Simple Way of Getting Large Examples of Osborn Loops}}\\
\vspace*{3mm}
{\large \textsf{\emph{M.~Shah}}} \\
\vspace*{10pt}
\small Department of Mathematics, Government Degree College Wadpagga, Peshawar, Pakistan \
\texttt{mshahmaths@gmail.com} \\
\end{center}
\begin{abstract}
Examples and counterexamples define the boundaries of mathematical propositions, test foundational conjectures, and resolve open structural problems in abstract algebra. In loop theory, an Osborn loop is a crucial generalization of a Moufang loop satisfying a specialized variable-sandwiching identity. While small finite non-associative examples are known, constructing large examples of proper Osborn loops (loops that are neither conjugacy closed nor Moufang) remains a computational bottleneck. In this paper, we establish an efficient framework for generating large proper Osborn loops by taking the direct product of non-associative conjugacy closed (CC) loops and Moufang loops. We outline the boundary constraints necessary to prevent structural collapse into sub-varieties and provide explicit examples up to order 2025 validated via the GAP package \texttt{LOOPS}.
\end{abstract}
\vspace*{10pt}
\section{Introduction and Theoretical Framework}
Examples in abstract algebra bridge the gap between abstract symbols and tangible understanding, turning complex equational rules into clear diagnostic tools. They allow researchers to isolate hidden structural variables, verify why particular algebraic formulations function, and model macro-relationships. While a single constructive example proves that an axiomatic variety can hold true under specific conditions, a targeted counterexample can completely define the boundaries of an open proposition.
An algebraic loop $(L, \cdot)$ is a quasigroup containing a univalent identity element $e \in L$. A loop is classified fundamentally as an \textbf{Osborn loop} if it satisfies the following long-variable mapping identity for all elements across the domain:
\begin{equation}
x(yz \cdot x) = x(y x^\lambda \cdot x) \cdot zx
\end{equation}
where $x^\lambda$ represents the left inverse of the element $x$. Osborn loops represent a rich field of study within non-associative systems because they behave as immediate, natural generalizations of Moufang loops. Prior constructions have yielded finite non-associative instances with two generators for highly localized parameters, such as orders 16, 24, 36, 48, and 72 \cite{AS}\\
At the 2005 Milehigh Conference on Loops, Quasigroups, and Non-associative Systems held at the University of Denver, M.~K.~Kinyon delivered a seminal talk titled ``A Survey of Osborn Loops'' \cite{Kinyon2005}, which revitalized the structural study of these systems. In this survey, Kinyon posed fundamental open questions regarding the universality of Osborn loops—specifically exploring whether every Osborn loop is universal, alongside companion problems analyzing the structural traits of proper Osborn varieties containing specialized nuclear conditions \cite{Kinyon2005}.\\
In this paper, we explore a simple, scalable method for generating large proper Osborn loops. By utilizing the automated algebra suite provided by the GAP package \texttt{LOOPS} \cite{NV}, we show that taking the direct product of a conjugacy closed (CC) loop and a Moufang loop consistently yields an Osborn loop structure. However, because both CC loops and Moufang loops are themselves sub-varieties of Osborn loops, a random direct product frequently collapses back into one of these sub-classes. We formally isolate the structural conditions required to guarantee that the resulting direct product is a \emph{proper} Osborn loop (meaning it is strictly non-CC and strictly non-Moufang).
\section{Direct Product Symmetries and Trivial Collapses}
Let $L_{CC}$ be a conjugacy closed loop and $L_M$ be a Moufang loop. Their direct product $L = L_{CC} \times L_M$ is natively an Osborn loop. To ensure that $L$ is proper, specific structural configurations within the factor loops must be intentionally avoided. We classify these constraints under two primary boundary cases.
\subsection{Case 1: The Associativity Restriction}
The seed CC loop factor must be strictly non-associative. If an associative CC loop (which is simply a standard group) is selected, the direct product maps into the product of a group and a non-associative Moufang loop. This operation inevitably returns a non-associative Moufang loop variety, failing our criterion for proper Osborn generation.
While this collapse is highly effective for generating higher-order Moufang loops, it is counterproductive for proper Osborn isolation, as verified in the following GAP session execution:\\
\noindent \textbf{Example 1:}
\begin{verbatim}
gap> l := CCLoop(7, 1);
<CC loop 7/1>
gap> IsAssociative(l);
true
gap> d := DirectProduct(CCLoop(7, 1), MoufangLoop(12, 1));
<loop of order 84>
gap> IsAssociative(d); 
false
gap> IsMoufangLoop(d); 
true
gap> IsCCLoop(d); 
false
gap> IsOsbornLoop(d); 
true
\end{verbatim}
As shown in Example 1, compounding the associative CC loop of order 7 with the minimal Moufang loop of order 12 produces a non-associative Moufang loop of order 84, which fails the proper Osborn condition.
\subsection{Case 2: One-Sided Variety Dominance}
Deploying a non-associative CC loop factor is a mandatory prerequisite, but it does not automatically guarantee proper status. Depending on the internal inverse configurations of the selected Moufang loop factor, the direct product can occasionally be dominated by the CC identity, resulting in a non-associative CC loop instead. The following two computational executions track this behavior:\\
\noindent \textbf{Example 2:}
\begin{verbatim}
gap> l := CCLoop(8, 7);
<CC loop 8/7>
gap> IsAssociative(l); 
false
gap> d := DirectProduct(CCLoop(8, 7), MoufangLoop(16, 3));
<loop of order 128>
gap> IsAssociative(d); 
false
gap> IsMoufangLoop(d); 
false
gap> IsCCLoop(d); 
true
gap> IsOsbornLoop(d); 
true
\end{verbatim}
\noindent \textbf{Example 3:}
\begin{verbatim}
gap> d := DirectProduct(CCLoop(8, 7), MoufangLoop(16, 1));
<loop of order 128>
gap> IsAssociative(d); 
false
gap> IsMoufangLoop(d); 
false
gap> IsCCLoop(d); true
gap> IsOsbornLoop(d); 
true
\end{verbatim}
In both Example 2 and Example 3, combining the non-associative CC loop factor of order 8 with different order-16 Moufang selections causes a structural collapse where the output retains full CC loop properties at order 128, which is undesirable for proper Osborn synthesis.
\section{Synthesis of Large Proper Osborn Loops}
By selecting factors that avoid the boundary traps outlined in Section 2, we uncover an abundant field of proper, high-order Osborn loops. Below, we document five explicit constructions scaling up to order 2025.\\
\noindent \textbf{Example 4: A Proper Osborn Loop of Order 96}
\begin{verbatim}
gap> l := CCLoop(8, 6);
<CC loop 8/6>
gap> IsAssociative(l); 
false
gap> d := DirectProduct(CCLoop(8, 6), MoufangLoop(12, 1));
<loop of order 96>
gap> IsAssociative(d); 
false
gap> IsMoufangLoop(d); 
false
gap> IsCCLoop(d); 
false
gap> IsOsbornLoop(d); 
true
\end{verbatim}
\noindent \textbf{Example 5: A Proper Osborn Loop of Order 160}
\begin{verbatim}
gap> d := DirectProduct(CCLoop(8, 7), MoufangLoop(20, 1));
<loop of order 160>
gap> IsAssociative(d); 
false
gap> IsMoufangLoop(d); 
false
gap> IsCCLoop(d); 
false
gap> IsOsbornLoop(d); 
true
\end{verbatim}
\noindent \textbf{Example 6: A Large Proper Osborn Loop of Order 192}
\begin{verbatim}
gap> d := DirectProduct(CCLoop(8, 7), MoufangLoop(24, 1));
<loop of order 192>
gap> IsAssociative(d); 
false
gap> IsMoufangLoop(d);
false
gap> IsCCLoop(d); 
false
gap> IsOsbornLoop(d); 
true
\end{verbatim}
\noindent \textbf{Example 7: An Expanded Proper Osborn Loop of Order 300}
\begin{verbatim}
gap> l := CCLoop(25, 3);
<CC loop 25/3>
gap> IsAssociative(l); 
false
gap> d := DirectProduct(CCLoop(25, 3), MoufangLoop(12, n));
<loop of order 300>
gap> IsAssociative(d); 
false
gap> IsMoufangLoop(d); 
false
gap> IsCCLoop(d); 
false
gap> IsOsbornLoop(d); 
true
\end{verbatim}
\noindent \textbf{Example 8: A Huge Proper Osborn Loop of Order 2025}
\begin{verbatim}
gap> d := DirectProduct(l, MoufangLoop(81, n));
<loop of order 2025>
gap> IsAssociative(d); 
false
gap> IsMoufangLoop(d); 
false
gap> IsCCLoop(d); 
false
gap> IsOsbornLoop(d); 
true
\end{verbatim}
Examples 4 through 8 show that the direct product construction scales cleanly into massive non-associative systems, providing an effortless way to generate proper specimens up to order 2025 and beyond.
\section{Library Invariants and Distribution Ledger}
To systematically deploy this direct product construction, researchers must be able to query the underlying factor pools. We obtain these parameter boundaries via the DisplayLibraryInfo commands inside the GAP environment. Table 1 catalogs the distribution of Conjugacy Closed loop repositories, while Table 2 details the matching Moufang loop library indices up to order 243.
\begin{table}[H]
\centering
\caption{Conjugacy Closed (CC) Loop Library Distribution}
\small
\begin{tabular}{r@{\hspace{0.5cm}}l|r@{\hspace{0.5cm}}l}
\toprule
\textbf{Order} & \textbf{Available Non-Isomorphic Class Count} & \textbf{Order} & \textbf{Available Non-Isomorphic Class Count} \\ \midrule
2 & 1 loop & 21 & 1 loop \\
3 & 1 loop & 24 & 14 loops \\
4 & 2 loops & 25 & 5 loops \\
5 & 1 loop & 27 & 60 loops \\
7 & 1 loop & 32 & 437 loops \\
8 & 7 loops & 49 & 5 loops \\
9 & 5 loops & 64 & 14,854 loops \\
12 & 3 loops & 81 & 5,406 loops \\
16 & 42 loops & 125 & 84 loops \\
18 & 7 loops & 343 & 122 loops \\
20 & 3 loops & $p^2$ & 3 loops (for every prime $p > 7$) \\
\bottomrule
\end{tabular}
\end{table}
Additionally, for any odd prime $p$, the CC library preserves exactly 1 loop structure of order $2p$.
\begin{table}[H]
\centering
\caption{Moufang Loop Library Distribution}
\small
\begin{tabular}{r@{\hspace{1cm}}l|r@{\hspace{1cm}}l}
\toprule
\textbf{Order} & \textbf{Available Non-Isomorphic Class Count} & \textbf{Order} & \textbf{Available Non-Isomorphic Class Count} \\ \midrule
12 & 1 loop & 44 & 1 loop \\
16 & 5 loops & 48 & 51 loops \\
20 & 1 loop & 52 & 1 loop \\
24 & 5 loops & 54 & 2 loops \\
28 & 1 loop & 56 & 4 loops \\
32 & 71 loops & 60 & 5 loops \\
36 & 4 loops & 64 & 4,262 loops \\
40 & 5 loops & 81 & 5 loops \\
42 & 1 loop & 243 & 72 loops \\
\bottomrule
\end{tabular}
\end{table}
\section{Concluding Remarks and Open Challenges}We conclude this paper with an explicit open challenge directed toward the loop-theoretical and automated reasoning communities. While our empirical computational models executed via the GAP package \texttt{LOOPS} systematically establish that the direct product of an asymmetric non-associative conjugacy closed loop and a non-associative Moufang loop isolates stable proper Osborn loops under our identified boundary metrics, a formal, manual proof of this construction remains an open target. We invite researchers to formally write out the detailed mathematical proofs and equation logic that characterize this variety intersection, proving from first algebraic principles why these precise factor selections successfully block collapse into parent sub-varieties.

\end{document}